ON THE REDUCTION PROCESS OF NUCCI-REDUCE ALGORITHM FOR
COMPUTING NONLOCAL SYMMETRIES OF DYNAMICAL SYSTEMS:
A case study of the Kepler and Kepler- related problems.

I.F.Arunaye
Department of Mathematics and Computer Science
University of the West Indies, Mona
Kingston 7, Jamaica

INTRODUCTION

In 1994 Krause used the Kepler problem as a vehicle for the introduction of his concept of the complete symmetry group of a differential equation or systems of differential equations. Essentially the complete symmetry group of a differential equation is the group of the algebra of the Lie symmetries required to specify completely the differential equation [2, 1]. In the case for the Kepler problem, the group comprises of nonlocal symmetries although it is not the case for every differential equations [1, 12, 5]. These nonlocal symmetries are not Noether symmetries (Variational symmetry). Since Krause (1994) the concern for nonlocal symmetry became prominent in the domain of symmetry analysis of dynamical systems [11, 5, 13, 8]. In 1996, Nucci unveil a profound algorithm for computing the nonlocal symmetries of the Kepler problem where she combined her reduce technique with the Lie point algorithm, and thereby invalidated Krause (1994) assertion that the complete symmetry group of the Kepler problem could not be ascertained by the Lie group analysis. Since then a substantial literature is devoted to topics in nonlocal symmetry analysis. In this regards the NUCCI-REDUCE algorithm software is acclaimed the most effective and efficient of all the techniques available in the literature to date [11, 2, 4,3]. However, attempting to use this REDUCE algorithm to compute the non-trivial problems with pencil and papers could be very complicated and complex. It suffices to say that one must develop ingenuity in the choice of change of variables and computational skills to reduce the nonlinear system to a system of two linear equations with one second order and $(n-1)$ first order (for an n-dimensional argument) equations which we now apply the Lie group algorithm to obtain the point symmetries in the final new variables. These are transformed back to the original variables to obtain both point and nonlocal symmetries of the original problem. The Nucci REDUCE algorithm had been used successfully to compute both types of symmetries for all other problems possessing conserved vectors of the Laplace-Runge-Lenz and Poincaré type, and even for Lagrangian systems and Hamiltonian equations [7, 12, 11, 6]. It is also true that the Nucci- reduce algorithm is successfully applicable to any system of n first-order equations, which are transformable into equivalent systems where at least one of the equations is of the second order and the admitted Lie symmetry algebra is no longer infinite dimensional and the nontrivial symmetries of the original system could be retrieved [9, 3, 14]. This paper is organized in the following manner; firstly it recalled the Nucci reduction algorithm with associated deficiencies in sections **1.0** and **2.0**, followed by the reduction technique we proposed which we applied to the classical Kepler and Kepler problem with drag force in **3.0** and more general form of the inverse law equation in **4.0.** We further apply this to other physical problems including the MICZ problem in sections **5.0-6.0** and obtained the symmetry group of the MICZ problem on the cone. Finally, it ends with concluding remarks.

**1.0   Nucci reduction process for the two dimensional Kepler problem**

$$\underline{\ddot{r}} = -\frac{\mu \underline{r}}{r^3}; \qquad |\underline{r}| = r \tag{1.1}$$

We denote by $(r, \theta)$ polar coordinates of the particle in the plane of motion.
The radial and transverse components of the equation of motion for (1.1)

$$\ddot{r} - r\dot{\theta}^2 = -\frac{\mu}{r^2} \tag{1.2}$$

$$r\ddot{\theta} + 2\dot{r}\dot{\theta} = 0$$

## Nucci algorithm

$$w_1 = r \quad ; \quad \dot{w}_1 = w_3$$

$$w_2 = \theta \quad ; \quad \dot{w}_2 = w_4$$

$$w_3 = \dot{r} \quad ; \quad \dot{w}_3 = w_1 w_4^2 - \frac{\mu}{w_1^2} \tag{1.3}$$

$$w_4 = \dot{\theta} \quad ; \quad \dot{w}_4 = \frac{2 w_3 w_4}{w_1}$$

since $\theta$ is ignorable coordinate in (1.1), select $w_2 = y$ as new independent variable. Then,

$$\frac{dw_1}{dy} = \frac{w_3}{w_4} \tag{1.4}$$

$$\frac{dw_3}{dy} = w_1 w_4 - \frac{\mu}{w_1^2 w_4} \tag{1.5}$$

$$\frac{dw_4}{dy} = -\frac{2 w_3}{w_1} \tag{1.6}$$

from (1.4) $w_3 = w_4 w_1'$ \hfill (1.7)

$$w_4 = -\frac{w_4 w_1'}{w_1} \tag{1.8}$$

from (1.5)

$$w_4 w_1'' - \frac{2 w_4 w_1'^2}{w_1} = w_1 w_4 - \frac{\mu}{w_1^2 w_4} \tag{1.9}$$

We note that the integrating factor $w_1^2$ makes (1.8) exact so,

$$(w_1^2 w_4)' = 0 \tag{1.10}$$

This prompts the choice of a new variable to be

$$u_2 = w_1^2 w_4 \tag{1.11}$$

Substituting (1.11) into (1.9)

$$u_2^2 \left(\frac{w_1''}{w_1^2} - \frac{2w_1'^2}{w_1^3}\right) = \frac{u_2^2}{w_1} - \mu \tag{1.120}$$

so, taking $u_1 = \mu - \dfrac{u_2^2}{w_1}$ \hfill (1.13)

Putting (1.13), 1.11) in (1.12) reduce the Kepler problem (1.1) to the system

$$u_1'' + u_1 = 0 \tag{1.14}$$

$$u_2' = 0$$

## 2.0 Nucci reduction Process for the kepler problem with drag

$$\ddot{\underline{r}} + \frac{\alpha \dot{r}}{r^2} + \frac{\mu r}{r^3} = 0 \tag{2.1}$$

The radial and transverse components of the motion are as follows

$$\ddot{r} - r\dot{\theta}^2 + \frac{\alpha \dot{r}}{r^2} + \frac{\mu}{r^2} = 0 \tag{2.2}$$

$$r\ddot{\theta} + 2\dot{r}\dot{\theta} + \frac{\alpha \dot{\theta}}{r} = 0 \tag{2.3}$$

The algorithm yields the following equations for the $w_i$:

$$w_1' = \frac{w_3}{w_4} \quad ; w_3' = w_1 w_4 - \frac{\alpha w_3}{w_1^2 w_4} - \frac{\mu}{w_1^2 w_4} \quad ; w_4' = -\frac{2w_3}{w_1} - \frac{\alpha w_4}{w_1^2}.$$

where $w_1' = w_{,y}, y = w_2 = \theta$.

Using the same procedure above to eliminate $w_3$, we obtain

$$w_1'' w_4 + w_1' w_4' = w_1 w_4 - \frac{\alpha w_1'}{w_1^2} - \frac{\alpha}{w_1^2 w_4} \tag{2.3}$$

$$w_4' = -\frac{2w_1'w_4}{w_1} - \frac{\alpha}{w_1^2} \tag{2.4}$$

$$(w_1^2 w_4)' + \alpha = 0 \quad \Rightarrow \quad w_1^2 w_4 = -(\alpha y + \beta) \tag{2.5}$$

Choose

$$u_2 = w_1^2 w_4 + \alpha y + \beta \tag{2.6}$$

Eliminating $w_4$ and $w_4'$ from (2.3), we obtain

$$\left(\frac{1}{w_1}\right)'' + \frac{1}{w_1} = \frac{\mu}{(\alpha y + \beta)^2} \tag{2.7}$$

So, we conveniently choose

$$u_1 = \frac{1}{w_1} + \mu \int^y \frac{\sin(y-s)}{(\alpha s + \beta)^2} ds \tag{2.8}$$

and equation (2.6) and (2.8) reduced to the system

$$u_1'' + u_1 = 0$$
$$u_2' = 0 \tag{2.9}$$

As in the Kepler problem, this system also reduced to a harmonic oscillator and a conservation law. This similarity between the two systems (1.1) and (2.1) has already been reported [13]. We note that where the equation of motion under consideration using Nucci reduction process has more than two independent variables, the reduction process become cumbersome [7]. The cases of the Kepler problem in 3-dimensions and the MICZ problem are simple acid test cases. Nucci algorithm involved a set of six first order equations which are further transformed into five equations (one second order and four first order) [11,12]. Then we required five new and a set of three variables from which a choice of ignorable variable is made. Further, we are required to eliminate from a set of four first order equations a certain variable to be able to obtain the z- component of the parabolic determining equation which variables provide new dependent variables that reduce the system to two second order equations and a conservation law. The expressions involved in these processes are not simple. It is an uphill task to derive the first derivatives of the variables involved with pencil and papers [7]. Consequently, the reduction process required high level of symbolic computer manipulations code that are computer time consuming.

The algorithm to find first integrals which are not necessarily "physical" conserved quantities, suffice it to say that there is need to adjust the algorithm of the reduction of order introduced by Nucci since all symmetries are found by Lie point method, and there exist natural variables with which to reduce the equation of motion to suitable system of linear equations and conservation laws [13, 6].

**3.0 Our proposed reduction process**

We construct an isomorphic transformation $\Phi(u) \cong F(r)$ such that the transformation is invariant, respect the constants of the motion, and preserved the orbit of the motion. Equation (1.2) is rewritten as

$$\ddot{r} - r\dot{\theta}^2 = -\mu r^{-2} \tag{3.1}$$

$$r^{-1}(r^2\dot{\theta})^{\cdot} = 0$$

$$r^2\dot{\theta} = L \equiv \text{angular momentum} \tag{3.2}$$

$$\phi(u) \equiv u = r^{-1} \ ; \ \dot{r} = -Lu_\theta \ ; \ \ddot{r} = -L^2 u^2 u_{\theta\theta} \tag{3.3}$$

$$-L^2 u^2 u_{\theta\theta} - u^3 L^2 = -\mu u^2 \tag{3.4}$$

$$L^2 u^2 (u_{\theta\theta} + u + -\mu_2 L^{-2}) = 0 \tag{3.5}$$

$$u_{\theta\theta} + u = \mu L^{-2} \tag{3.6}$$

Defining variables $u_1 = u - \mu L^{-2}$, $u_2 = L$, then we obtain the system of equations

$$u_{1,\theta\theta} + u_1 = 0 \tag{3.7}$$

$$u_{2,\theta} = 0$$

This is a system of two equations with one second order (harmonic oscillator) and a conservation law that we can be easily solved by Lie point symmetry algorithm to have nine point symmetries. Just as well known, we can translate back to original variables to obtain the symmetries, some point, others nonlocals. That is,

$\bar{\theta} = \theta + \lambda \xi(\theta, u_i)$, $\bar{u} = u + \lambda \eta_i(\theta, u_i)$; with symmetry generator

$V = \xi \partial_\theta + \eta_i \partial_i$ where $\partial_i = \partial/\partial u_i$.

With second prolongation given by

$$V^{[2]} = V + \phi_i^1 \partial_{u_{i,\theta}} + \phi_i^2 \partial_{u_{i,\theta\theta}}$$

where $\phi_i^1, \phi_i^2$ are the first and second prolongation functions, produced the symmetry equations that are solved using the method of superimpositions to have eight point symmetries while $u_{2,\theta}$ account for the time translation symmetry just as in the Nucci reduce algorithm. Using the same procedure for Kepler problem with drag force, we have the force components as:

$$\ddot{r} - r\dot{\theta}^2 + \frac{\alpha \dot{r}}{r^2} = -\frac{\mu}{r^2}$$

$$r\ddot{\theta} + 2\dot{r}\dot{\theta} + \frac{\alpha \dot{\theta}}{r} = 0 \implies \dot{L} + \alpha\dot{\theta} = 0$$

i.e. $L + \alpha\theta$ is a constant.

Using the substitution $u = r^{-1}$ as in the Kepler problem we note tha

$$\dot{r} = -u^{-2}\dot{\theta}u_\theta = -Lu_\theta$$

$$\ddot{r} = -L^2 u^2 \dot{\theta} u_{\theta\theta} - \dot{L}u_\theta. \tag{3.8}$$

Hence the radial equation becomes

$$-L^2 u^2 u_{\theta\theta} - \dot{L}u_\theta - \alpha u^2 L u_\theta + u^3 L^2 = -\mu u^2$$

i.e. $\quad u_{\theta\theta} + u = \mu L^{-2}. \tag{3.9}$

A particular solution to this equation is

$$\upsilon = \mu \int^\theta \sin(\theta - \eta)(L_\circ - \alpha\eta)^{-2} d\eta \tag{3.10}$$

where $L_\circ$ is the constant $L + \alpha\theta$. Defining variables $u_1 = u - \upsilon$, $u_2 = L + \alpha\theta$, then as in the Kepler problem, we obtain the system

$$u_{1,\theta\theta} + u_1 = 0 \tag{3.11}$$

$$u_{2,\theta} = 0.$$

**4.0 A more general equation of motion in 2-dimensions**

(Plane polar coordinate system)

The general equation

$$\ddot{\mathbf{r}} + \mu r^\alpha \mathbf{r} = 0 \tag{4.1}$$

has two componentes of the motion in the form

$$\ddot{r} - r\dot{\theta}^2 = -\mu r^{\alpha+1} \tag{4.2}$$

$$r\ddot{\theta} + 2\dot{r}\dot{\theta} = 0 \tag{4.3}$$

from (4.3), $(r^2\dot{\theta})^\cdot = 0 \Rightarrow r^2\dot{\theta} = L \equiv$ constant

from (4.2), we have by our usual reduction

$$-L^2 u^2 u_{\theta\theta} - u^3 L^2 = -\mu u^{-(\alpha+1)}$$

$$L^2u^2[u_{\theta\theta}+u-\mu L^{-2}u^{-(\alpha+3)}]=0$$

$$u_{\theta\theta}+u-\mu L^{-2}u^{-(\alpha+3)}=0 \qquad (4.4)$$

take $\alpha+3=k$

The case $k=0$, then $\alpha=-3$, is the Kepler problem of (3.7).

For $k=-1$, then $\alpha=-4$, we have

$$u_{\theta\theta}+(1-\mu L^{-2})u=0$$

and the system (4.3) becomes

$$u_{1,\theta\theta}+(1-\mu L^{-2})u_1=0 \qquad (4.5)$$

$$u_{2,\theta}=0$$

This is the reduced form of the inverse cube law problem.

## 5.0 Generalization of the Kepler problem with drag on the cone

The equation of motion is of the form

$$\ddot{\mathbf{r}}-(\frac{\dot{g}}{2g}+\frac{3\dot{r}}{2r})\dot{\mathbf{r}}+\mu g\mathbf{r}=0 \qquad (5.1)$$

The radial and transverse components of equation of motion in the plane of motion are

$$\ddot{r}-r\dot{\theta}^2=-(\frac{\dot{g}}{2g}+\frac{3\dot{r}}{2r})\dot{r}-\mu gr \qquad (5.2)$$

$$r^{-1}(r^2\dot{\theta})^{\cdot}=(\frac{\dot{g}}{2g}+\frac{3\dot{r}}{2r})r\dot{\theta} \qquad (5.3)$$

Equation (5.3) implies that

$$\dot{L}=\left(\frac{\dot{g}}{2g}+\frac{3\dot{r}}{2r}\right)L=\tfrac{1}{2}(\ln gr^3)^{\cdot}L \qquad (5.4)$$

where $L=r^2\dot{\theta}$. Thus

$$L=A(gr^3)^{1/2} \qquad (5.5)$$

where $A$ is a constant.

Using the formulae in (3.8) for $\ddot{r}, \dot{r}$ etc, the radial equation (5.2) becomes

$$-L^2 u^2 (u_{\theta\theta} + u) + \mu g r = 0$$

i.e. $\quad A^2(u_{\theta\theta} + u) - \mu = 0$

$$u_{\theta\theta} + u = \mu A^{-2}. \tag{5.6}$$

Defining the two variables $u_1 = u - \mu A^{-2} = u - \mu u_2^{-2}$, $u_2 = L(gr^3)^{-\frac{1}{2}}$, we obtain the linear system

$$u_{1,\theta\theta} + u_1 = 0 \tag{5.7}$$

$$u_{2,\theta} = 0.$$

## 6.0 The MICZ problem on the cone of motion

The equation of motion is given by

$$\ddot{\mathbf{r}} + \frac{\lambda \mathbf{L}}{r^3} + (\frac{\mu}{r^3} + \frac{2\nu}{r^4})\mathbf{r} = 0; \quad \text{for} \quad -\lambda^2 = 2\nu \quad \text{we have}$$

$$\ddot{\mathbf{r}} + \frac{\lambda \mathbf{L}}{r^3} + (\frac{\mu}{r^2} - \frac{\lambda^2}{r^3})\hat{\mathbf{r}} = 0 \tag{6.1}$$

$$\ddot{r} - r\dot{\phi}^2 \sin^2\theta = -\frac{\mu}{r^2} + \frac{\lambda^2}{r^3}$$

$$\ddot{r} - L^2 r^{-3} = -\frac{\mu}{r^2} + \frac{\lambda^2}{r^3}, \quad \text{where } L = r^2 \sin\theta\dot{\phi} \text{ and } (\sin\theta)^{-1}L = u^{-2}\dot{\phi}$$

$$\ddot{r} - (\lambda^2 + L^2)r^{-3} = -\mu r^{-2} \tag{6.2}$$

setting $u = r^{-1}$, $\dot{r} = -u^{-2}\dot{u} = -u^{-2}u_\phi \dot{\phi}$

$$\ddot{r} = -L(\sin\theta)^{-1} u_{\phi\phi} \dot{\phi} = -L^2(\sin\phi)^{-2} u_{\phi\phi} u^2$$

(6.2) implies

$$-L^2(\sin\phi)^{-2} u_{\phi\phi} u^2 - (\lambda^2 + L^2)u^3 = -\mu u^2$$

$$u_{\phi\phi} + S^2(1+\lambda^2 L^{-2})u = \mu S^2 L^{-2} \tag{6.3}$$

where $S = \sin\theta$

But Poincaré vector gives

$$P^2 = L^2 + \lambda^2, \quad \mathbf{P}.\mathbf{r} = -\lambda r, \quad P\cos\theta = -\lambda, \quad \cos\theta = \frac{-\lambda}{(L^2+\lambda^2)^{\frac{1}{2}}} \text{ and }$$

$\sin\theta = \dfrac{L}{(L^2+\lambda^2)^{\frac{1}{2}}}$ . So, (6.3) reduced to

$$u_{\phi\phi} + u = \frac{\mu}{L^2 + \lambda^2} \tag{6.4}$$

and the MICZ problem reduces to

$$u_{1,\phi\phi} + u_1 = 0 \tag{6.5}$$

$$u_{2,\phi} = 0$$

where $u_1 = u - \dfrac{\mu}{L^2+\lambda^2}, u_2 = L$

This possess nine point symmetries which when transformed to the original variables becomes nonlocal ones.

In the case where $2v \neq -\lambda^2$ and is arbitrary, the MICZ reduces to

$$u_{\phi\phi} + \Omega^2 u = K = \frac{\mu}{L^2 + \lambda^2} \tag{6.6}$$

where $\Omega^2 = S^2(L^2 - 2v)$

If we set $\Omega\phi = x$ then (6.6) becomes

$$u_{xx} + u = K = \Omega^{-2}\frac{\mu}{L^2+\lambda^2} \tag{6.7}$$

The pair of equations corresponding to (6.5) is

$$u_{1,xx} + u_1 = 0$$

$u_{2,x} = 0$

where $u_1$ is now $u - \Omega^{-2} \dfrac{\mu}{L^2 + \lambda^2}$.

This system also possesses nine symmetry generators. However the explicit expressions for the generators are much more complicated and we do not proceed further in this case.

### 6.1 The symmetry group of the MICZ on the cone is obtained by the following transformation:

A $(t, r, \phi)$ symmetry becomes a $(\phi, u_1, u_2)$ symmetry according to following transformation

$$V = \xi^t \partial_t + \xi^r \partial_r + \xi^\phi \partial_\phi \to W = \sigma \partial_\phi + \Omega \partial_{u_1} + \Sigma \partial_{u_2} \qquad (6.8)$$

We note that $u_2 = r^2 \dot{\phi}$, $u_1 = \dfrac{1}{r} - \dfrac{\mu}{u_2^2 + \lambda^2}$ and by considering the relations

$W(\phi) = V(\phi)$, $W(u_1) = V(r^{-1} - \dfrac{\mu}{L^2 + \lambda^2})$ etc, one obtains the following

$$\Sigma = 2\xi^r r\dot{\phi} + r^2(\dot{\sigma} - \dot{\phi}\dot{\xi}^r)$$

$$\Omega = -\dfrac{\xi^r}{r^2} + \dfrac{2\mu u_2}{(u_2^2 + \lambda^2)^2} \Sigma. \qquad (6.9)$$

The representations of the complete symmetry group of the Kepler problem abound in the literature. If in the symmetry representation of the Kepler problem [8], we replace $L^2$ by $L^2 + \lambda^2$ we obtain the complete symmetry group representation of the MICZ problem from (6.8) and (6.9) on its cone of motion as following except $\Lambda_1$ which is a completely different from the scaling symmetry from that of the Kepler problem:

$$\Lambda_1 = (-3t + \int \dfrac{4\mu\lambda^2 r}{L_1^2} dt)\partial_t + (-2r + \dfrac{2\mu\lambda^2}{L_1^2} r^2)\partial_r$$

$$\Lambda_2 = \partial_\phi$$

$$\Lambda_3 = 2[\mu\int r dt - L_1^2 t]\partial_t + (\mu r - L_1^2)\partial_r$$

$$\Lambda_{4\pm} = 2[\int r e^{\pm i\phi} dt]\partial_t + r^2 e^{\pm i\phi}\partial_r$$

$$\Lambda_{6\pm} = 2[\int (\mu r + 3L_1^2)e^{\pm i\phi} dt]\partial_t + r(\mu r + 3L_1^2)e^{\pm 2i\phi}\partial_r + L_1^2 e^{\pm 2i\phi}\partial_\phi$$

$$\Lambda_{8\pm} = 2[\int \{2\dot{r}L_1^3 \pm ir(\mu - r^{-1}L_1)(\mu + r^{-1}L_1)\}e^{\pm i\phi} dt]\partial_t + r[2\dot{r}L_1^3 \pm ir(\mu - r^{-1}L_1)(\mu + r^{-1}L_1)]\partial_r$$

$$+ L_1^2(\mu - r^{-1}L_1)e^{\pm i\phi}\partial_\phi$$

where $L_1^2 = L^2 + \lambda^2$ and note that for $\lambda = 0$, we have the Kepler problem.

In the Kepler case, the two nonlocal symmetries to be included to form the complete symmetry group are of type $\Lambda_{4\pm}$. We have not seen where this has been reported.

**Concluding remarks**

It is common knowledge that all point symmetries can be determined by Lie group algorithm; there is also the potential for a loss of symmetries in the reduction process just as there is hope for an increase in the total number of symmetries [13]. The algorithm for translating back to symmetries in the original variables of the unreduced system is not ambiguous as it where [13]. Sine-qua-non, the concerns of the reduction process ought to be more centered on reducing the system to a set of equations which one can easily apply Lie algorithm to, rather than involving extraneous ambiguities such as determining first integrals and complexities of the choice of varying variables that involved complex computer manipulations and computer time wastage.

In this paper we have been able to establish with examples a simple algebraic process for reducing systems to equivalent systems which Nucci algorithm could obtain that are easily solved for their symmetries by Lie group analysis without ambiguities. We report here that we have used this algebraic process to reduce other nonlinear systems to convenient forms which we successfully applied the Lie algorithm to and we obtained the Lie point symmetries and subsequently obtained the nonlocal symmetries in the original systems, both in two and three dimensional coordinate systems. To mention but a few, the Kepler-Ermakov, Ermakov-Toy, and any physical problem possessing the Poincaré vector are examples that we have reduced with this algebraic process.

The analysis of their symmetry algebras of equations of motion on the cone is a subject for discussion in another paper.

**Acknowledgement**

I hereby thank the University of the West Indies, Mona Campus, Jamaica ; for her support during this research and for funding my trip to Kiev-Ukraine to present it at the Seventh international Conference "Symmetry in Nonlinear Mathematical Physics" in June, 2007. I also thank my Department for their support in the course of this work.

**References**

1. Andriopoulos K., Leach P.G.L and Flessas G.P. Complete symmetry Groups of ordinary Differential Equations and Their Integrals: Some Basic Considerations. J. Math. Anal. and Appl. 262, 256-273. 2001.

2. Andriopoulos K and Leach P.G.L. The Economy of Complete symmetry Groups for linear Higher Dimensional systems. J. Nonlinear Math. Phys. 9, Sup.2 ,10-23, 2002

3. Edwards M and Nucci M.C. Application of Lie group analysis to a core group model for sexually transmitted diseases. J. Nonlinear Maths. Phys. 13, N2, 211-230, 2006.

4. Gandarias M.L. Medina E. and Muriel C. New symmetry Reductions for some ordinary Differential Equations. J.Nonlinear Math. Phys. 9, Sup. 1, 47-59, 2002.

5. Krause J. On the complete symmetry group of the classical Kepler system. J. Math. Phys. 35, N11, 5735-5748,1994.

6. Leach P.G.L, Andriopoulos K. and Nucci M.C. The Ermanno-Bernoulli constants and representations of the complete symmetry group of the Kepler problem. J. Math. Phys. 44,N9, 4090-4106, 2003.

7. Leach P.G.L and Nucci M.C. Reduction of the classical MICZ problem to a two-dimensional linear isotropic harmonic oscillator. J.Math. Phys. 45, N9, 3590-3604, 2004.

8. Leach P.G.L and Flessas G.P. Generalizations of the Laplace-Runge-Lenz vector. J. Nonlinear Math. Phys. 10,N3, 340-423, 2003.

9. Marcelli M and Nucci M.C. Lie Point symmetries and first Integrals.: The Kowalevski top. J. Math. Phys. 44,N5, 2111-2131, 2003.

10. Nucci M.C. "Interactive REDUCE programs for calculating Lie point, nonclassical,Lie-Bäcklund, and approximate symmetries of differential equations: manual and floppy disk" CRC Handbook of Lie Group Analysis of Differential Equations edited by N.H.Ibragimov (CRC press, Boca Raton,1996), vol. III, pp415-481.

11. Nucci M.C. The complete Kepler group can be derived by Lie group analysis. J. Math. Phys. 37, N4, 1772-1775, 1996.

12. Nucci M.C. and Leach P.G.L. The Determination of Nonlocal symmetries by the Technique of Reduction of Order. J. Math. Anal. and Appl. 251, 871-884, 2000.

13 Nucci M.C and Leach P.G.L The harmony in the Kepler and related problems. J. Math. Phys. 42, N2, 746-764, 2001.

14. Nucci M.C. Lorenz integrable system Moves à la Poinsot. J. Math. Phys. 44, N9, 4109-4118,2003.